\documentclass[a4paper]{article}

\usepackage[utf8]{inputenc}
\usepackage{lmodern}
\usepackage{amssymb,amsfonts,amsthm,amsmath,bbm,mathrsfs,aliascnt,mathtools}
\usepackage[english]{babel}
\usepackage[colorlinks]{hyperref}

\theoremstyle{plain}
\newtheorem{theorem}{Theorem}
\newtheorem{corollary}[theorem]{Corollary}
\newtheorem{lemma}[theorem]{Lemma}

\newtheorem{remark}[theorem]{Remark}

\newcommand{\R}{\mathbb{R}}

\newcommand{\calS}{\mathcal{S}}

\newcommand{\calN}{\mathcal{N}_S}

\newcommand{\ind}[1]{\mathbf{1}_{\left\{#1\right\}}}

\renewcommand{\bar}[1]{\overline{#1}}

\newcommand{\e}{\mathrm{e}}
\newcommand{\dd}{\mathrm{d}}

\DeclareMathOperator{\E}{\mathbb{E}}

\renewcommand{\P}{\mathbb{P}}

\title{Limits of P\'olya urns with innovations}
\author{Jean Bertoin\thanks{Institute of Mathematics, University of Zurich, Switzerland.} }
\date{ }

\begin{document}

\maketitle

\begin{abstract} We consider a version of the classical P\'olya urn scheme which incorporates innovations. The space $S$ of colors is an arbitrary measurable set. 
After each sampling of a ball in the urn, one returns $C$ balls of the same color and additional balls of different colors given by some finite point process $\xi$ on $S$, where the distribution $P_s$ of the pair $(C,\xi)$ depends on the sampled color $s$. We suppose that the average number of copies $E_s(C)$ is the same for all $s\in S$, and that the intensity measures of innovations have the form $E_s(\xi)=a(s)\mu$ for some finite measure  $\mu$ and a  modulation function $a$ on $S$ that is bounded away from $0$ and $\infty$. We then show that  the empirical distribution of the colors in the urn converges to the normalized intensity  $\bar \mu$.  In turn, different regimes for the fluctuations are observed, depending on whether $E_s(C)$ is larger or smaller than $\mu(a)$. 
\end{abstract}

\noindent \emph{\textbf{Keywords:} P\'olya urn, innovation, empirical distribution, martingale central limit theorem.} 

\medskip

\noindent \emph{\textbf{AMS subject classifications:}} 60F17, 60G44, 60J85, 62G30.

\section{Introduction}

At  each step of an original P\'olya urn scheme,  a ball is drawn uniformly at random from an  urn, independently of the preceding steps. One observes its color; then  the ball is returned in the urn together with an additional ball of the same color.  More generally, one can consider random replacement schemes, for which one rather returns random numbers of balls of different colors according to some fixed distribution that depends only on the color of the sampled ball.
The quantities of interest are the proportions of balls with given colors after a large number of steps. 
Most often, it is assumed that the set of all possible colors  is finite, with the notable exception of certain recent contributions which will be discussed at the end of this Introduction. 
General  P\'olya urn schemes are used as models in a variety of fields including Computer Sciences, Biology, Social Sciences, etc.; see e.g. \cite{Mahmoud} and references therein.

We are interested in a variation of the P\'olya urn scheme which incorporates innovations, in the sense that at each step, balls with new colors that have never been used before can be returned in the urn, and the space of colors is an arbitrary measurable space. Such an urn scheme with innovation has been first introduced\footnote{Simon model bears close connections to the earlier work of Udny Yule on the frequency distribution of the sizes of biological genera, and in turn has been  rediscovered and applied in a variety of  frameworks, see e.g. the discussion in \cite{Willis}. } in 1955 
by  Herbert Simon \cite{Simon} in terms of a simple model
 for producing a text  by random iterations.  One starts with a first word, say $w_1$, and when one has already written a text of length $k\geq 1$, say $w_1 \ldots w_k$,  one appends the next word  $w_{k+1}$
 which is either, with fixed probability $p$, a repetition of $w_{u(k)}$, where $u(k)$ has the uniform distribution on $\{1, \ldots, k\}$, or, with complementary probability $1-p$,   a new word different from all the preceding. Further situations where innovations or novelties are combined with urn schemes have been considered by  \cite{Hoppe, Triaetal, Iacopini}, amongst others. Needless to recall, the association of  stochastic reinforcement and innovation is an important aspect of machine learning on which there exists of course a huge literature.

 In a recent work \cite{BeKest} that was partly motivated by the study of certain step reinforced random walks, we analyzed the effects of linear reinforcement  on the sequence of empirical distributions in Simon model.  It was observed there that when one takes for the different words independent uniform random variables on $[0,1]$, the Glivenko-Cantelli theorem still holds, that is there is the almost-sure convergence of the empirical distribution functions of the reinforced sequence towards the identity function on $[0,1]$, but the Donsker theorem is affected. Indeed,  the sequence of empirical processes converges in law to a Brownian bridge only up to a constant factor when the repetition probability $p<1/2$; a further rescaling is needed  when  $p>1/2$ and the limit is then a bridge with exchangeable increments and discontinuous paths. 
 
 Here, we consider more generally  a Borel space $S$ as the set of possible colors of balls. 
We write $\mathcal{L}^{\infty}$ for the space of bounded measurable functions $f: S\to \R$, and use often use the short notation
 $$\nu(f)\coloneqq \int_{ S}f(s)\nu(\dd s)$$
 for a finite measure $\nu$ on $S$ and $f\in \mathcal{L}^{\infty}$.
In particular,  the total mass of $\nu$ is denoted by $\nu(\mathbf{1})$, where $\mathbf{1}$ stands for the function identical to $1$ on $S$.
Whenever $\nu \not\equiv 0$, we write 
$$\bar \nu ={\nu(\mathbf{1})}^{-1} \nu$$ 
for the normalized probability measure associated to $\nu$.

A typical replacement consists of a pair $(C, \xi)$, where $C$ is a random variable with values in $\{-1,0,1,2, \ldots\}$ which represents the number of copies of the sampled ball which are returned in the urn (the case $C=-1$ accounts for the situation where the sampled ball is removed from the urn), and $ \xi$  a point process on $ S$ which represents  the random family of new balls which are simultaneously added.  
 The dynamics are hence fully encoded by the kernel  of  laws $(P_s)_{s\in S}$  that specify the distributions of the  pair  $(C, \xi)$ as a function of the color $s$ of the sampled ball. 
 We implicitly assume that the probability kernel $s\mapsto P_s$ is measurable, and to avoid discussing situations where the urn could be emptied, we always assume henceforth  that 
 for all $s\in S$, 
 \begin{equation}\label{E:1}
 P_s\left( C+\xi(\mathbf{1})\geq 0\right) = 1,
 \end{equation}
   where, by the preceding notation, $\xi(\mathbf{1})$ accounts for the  number of innovative balls.

  More precisely,  at  each step, independently of the previous steps, we sample a ball uniformly at random from the urn and observe its color, say $s$, as well as a pair  $(C,\xi)$ with law $P_s$. 
Given $C=k$ and $ \xi=\sum_{i=1}^j \delta_{s_i}$, with  $\delta_s$ denoting the Dirac measure at  $s$, we then return that ball in the urn together with $k$ balls with the same color $s$ and $j$ further balls with colors $s_1, \ldots, s_j$ (it is not required that the $s_i$ are all distinct).
We represent the composition of the urn after $n$ steps by a counting measure $U_n$ on $ S$ such that the atoms of $U_n$ are the colors of the balls present in the urn after $n$ steps, and their multiplicities the number of balls of these colors. Then for $f\in \mathcal L^{\infty}$,  $U_n(f)$ is the sum of the $f(s)$ for the  colors $s\in S$ present in the urn after $n$ steps and repeated according to their multiplicities.

We now introduce the three crucial assumptions of our study.
First, we suppose that the average number of copies that are returned at a typical step does not depend on the sampled color, viz.
\begin{equation}\label{E:2}
\text{the function $s\mapsto E_s(C)$ is constant on $S$,}
\end{equation}
where the notation $E_s$ is used for the mathematical expectation under $P_s$. 
Second, we suppose that the intensity measures of innovations are all proportional to some fixed  measure $\mu$ on the space of colors.  Specifically, we assume that  there exists a  measurable function $a$ on $S$ that is bounded away from $0$ and from $\infty$,
viz. with 
$$0< \inf_{s\in S} a(s)  \leq \sup_{s\in S}a(s)\coloneqq \|a\|_{\infty} < \infty,$$
such that  for every $s\in S$, 
the intensity measure of the point processes $\xi$ under $P_s$ is given by
\begin{equation} \label{E:3}
E_s(\xi(f))= a(s) \mu(f), \qquad \text{for all }f\in \mathcal L^{\infty}.
\end{equation}
We should think of $a$ as a factor which modulates the intensity of innovations as a function of the color of the sampled ball.  Last, we shall also need an assumption
of uniform upper and lower boundedness for 
the total number of balls which are returned at each step, namely 
that
\begin{equation} \label{E:4}
0< \inf_{s\in S} P_s(C+\xi(\mathbf{1})\geq 1) \quad\text{and} \quad 
\sup_{s\in S} E_s\left(\left( C+\xi(\mathbf{1})\right)^2\right) < \infty.
\end{equation}

 We then introduce two fundamental constants,
 \begin{equation}\label{E:5}
 \lambda_1\coloneqq E_s(C) + \mu(a)
 \quad\text{and}\quad \lambda_2\coloneqq E_s(C).
\end{equation}
As the notation suggests, these two quantities shall be viewed as eigenvalues. More precisely, 
consider 
 the operator 
 $\mathcal R: \mathcal{L}^{\infty}\to \mathcal{L}^{\infty}$ defined by
\begin{equation} \label{E:6}
\mathcal Rf(s) \coloneqq E_s(C f(s) + \xi(f))= \lambda_2 f(s) + a(s) \mu(f),
\end{equation}
which describes averaged replacements as a function of the color of the sampled ball. 
Then $\lambda_1$ and $\lambda_2$ arise as the two eigenvalues of $\mathcal R$; specifically
\begin{equation} \label{E:7}
\mathcal R a=\lambda_1 a \quad \text{and} \quad \mu(\mathcal Rf) = \lambda_1 \mu(f)  \quad \text{for all }f\in \mathcal L^{\infty},
\end{equation}
and also
\begin{equation} \label{E:8}
\mathcal R f=\lambda_2 f\quad \text{ whenever }\mu(f)=0.
\end{equation}
Observe from \eqref{E:1} that
$$\lambda_1 =\int_{S} E_s(C+\xi(\mathbf{1})) \bar \mu(\dd s) >0,
$$
and also that $\lambda_1> \lambda_2$,
since $\mu(a)>0$. 

When  the set of colors $S$ is finite, say $S=\{1, \ldots , d\}$, the present model with innovation merely rephrases the classical setting of urns with random replacements schemes and the average replacement operator $\mathcal R$ is given by a $d \times d$ matrix. The requirements  \eqref{E:2} and \eqref{E:3} entail that the latter has the form $\mathcal R= \lambda_2 \mathrm{Id}+\mathcal{R}_1$, with  $\mathcal{R}_1$ a matrix of rank $1$.
In this situation, the spectral properties of  $\mathcal R$ observed above are immediate, and 
the results in this work should not come much as a surprise. Nonetheless dealing with a general space of colors
cannot be reduced  (e.g. by approximations) to the case when $S$ is finite. Although our analysis will of course borrow the same guiding lines as for the classical setting, in particular spectral analysis of the mean replacement operator and martingale limit theorems,  
new ideas are also needed to resolve several issues.

An important result for urn schemes with finitely many colors is that the first order asymptotic of the contain of the urn as the number of steps goes to infinity is determined by the largest eigenvalue of the mean replacement matrix and its eigenvectors. The same  feature holds in the present setting; it implies the almost-sure convergence of the empirical distribution of colors to the normalised  intensity of  innovation.

\begin{theorem}\label{T:LLN} Assume \eqref{E:2}, \eqref{E:3}, and \eqref{E:4}.
For every $f\in \mathcal{L}^{\infty}$, we have
$$ \lim_{n\to \infty} n^{-1} U_n(f) = \lambda_1 \bar \mu(f)\qquad  \text{a.s.},$$
 where $\bar \mu \coloneqq \mu(\mathbf{1})^{-1} \mu$ is the normalised  intensity of  innovation. 
 As a consequence, we have also in terms of  the empirical distribution of colors  that 
\begin{equation} \label{E:9} \lim_{n\to \infty} \bar U_n(f) = \bar \mu(f)\qquad  \text{a.s.} 
\end{equation}
 \end{theorem}

 We next turn our attention to the fluctuations of the empirical distributions; in order to state our main result, we need to introduce first a few more objects. For any $f \in \mathcal L^{\infty}$, we set
 \begin{equation} \label{E:10}
 \sigma^2(f) \coloneqq \int_{S} E_s(\bar \mu(f^2) C^2 + \xi(f)^2) \bar \mu(\dd s),
 \end{equation}
 which defines a positive semidefinite quadratic form on $\mathcal L^{\infty}$. 
 We denote by $G=(G(f): f\in \mathcal L^{\infty})$ the associated Gaussian process, i.e. $G$ is  centered with covariance
 \begin{equation} \label{E:11}
 \E(G(f)G(g)) =  \int_{S} E_s\left( \bar \mu(fg) C^2 + \xi(f)\xi(g) \right)) \bar \mu(\dd s).
 \end{equation}
  In general, the variables $G(f)$ and $G(g)$ are not always independent when $f$ and $g$ are orthogonal in $\mathcal L^2(\mu)$,  due the component $  \E(\xi(f)\xi(g))$ in the covariance. 
  We also introduce the Gaussian bridge $G^{\mathrm{(br)}}=(G^{\mathrm{(br)}}(f): f\in \mathcal L^{\infty})$, where
  $$G^{\mathrm{(br)}}(f)\coloneqq G(f) - \bar \mu(f) G(\mathbf{1}).$$
  Plainly, the Gaussian bridge is translation invariant in the sense that $G^{\mathrm{(br)}}(f)=G^{\mathrm{(br)}}(f+c)$ for any $c\in \R$, and coincides with $G$
  on the hyperplane of  functions $ f\in \mathcal L^{\infty}$ with $\mu(f)=0$. 
 
It is well-known for classical urn schemes with finitely many colors, that the fluctuations of the empirical distributions of colors in the urn
depend crucially of whether the largest eigenvalue of the mean replacement matrix is larger or smaller than twice the real part of the second largest eigenvalue.
 This incites us to introduce
   \begin{equation} \label{E:12}
 \rho\coloneqq \frac{\lambda_2}{\lambda_1}\in(-\infty,1).
 \end{equation}
In particular, $\rho>1/2$ if and only if $E_s(C)>\mu(a)$, which we interpret as reinforcement being stronger than innovation.

 \begin{theorem} \label{T1} Assume \eqref{E:2}, \eqref{E:3}, and \eqref{E:4}. Then
 the following limits hold:
 
  \begin{enumerate}
 \item[(i)] If $\rho>1/2$, then we have for every  $f\in \mathcal L^{\infty}$:
 $$\lim_{n\to \infty}  n^{1-\rho} \left(\bar U_n(f)-\bar \mu(f)\right) = L_f \qquad \text{a.s.,}$$
 where $L_f $ is some random variable which is non degenerate, except when $\mu\left(|f|\right)=0$.  

 \item[(ii)] If $\rho<1/2$,  and further a slightly stronger version of \eqref{E:4},
  \begin{equation} \label{E:14}
\sup_{s\in S} E_s\left(\left( C+\xi(\mathbf{1})\right)^{\alpha}\right) < \infty \qquad \text{for some } \alpha>2, 
\end{equation}
holds, then the sequence of processes 
 $$ \left( \lambda_1  \sqrt{ (1-2\rho)n } \left(\bar U_n(f)-\bar \mu(f)\right): f\in \mathcal L^{\infty}\right), \qquad n\geq 1$$
converges in the sense of finite dimensional distributions to the Gaussian bridge $G^{\mathrm{(br)}}$.

 \item[(iii)] If $\rho=1/2$ and  \eqref{E:14} holds,  then the sequence
  $$ \left(\lambda_1 \sqrt{  \frac{n}{\log n}} \left(\bar U_n(f)-\bar \mu(f)\right): f\in \mathcal L^{\infty}\right), \qquad n\geq 1$$
converges in the sense of finite dimensional distributions to the Gaussian bridge $G^{\mathrm{(br)}}$.
 \end{enumerate}
  \end{theorem}
  
  \begin{remark} {\normalfont  \begin{enumerate}
  \item[(i)] Theorem \ref{T1}(ii-iii)  should be viewed as weaker analogs of \linebreak Donsker's limit theorem for empirical distributions. It would be interesting to strengthen  
convergence in the sense of finite dimensional distributions here to a uniform central limit type theorem on appropriate classes of functions as in Dudley \cite{Dudley}.

\item[(ii)] For instance, Simon's model with words given by uniform random variables discussed at the beginning of this introduction corresponds to  the case where $S=[0,1]$, and $P_s$ is the law of a pair $(C,\xi)$ where $C$ is a Bernoulli variable with parameter $p$ and  $\xi =(1-C)\delta_V$, with $V$  a uniform variable on $[0,1]$ independent of $C$. hence $P_s$ does not depend on $s$, we may take 
 $\mu(\dd s) = (1-p)\dd s$, $a \equiv 1$, $\lambda_1=1$ and $\lambda_2=p$. One finds
 $\sigma^2(f) = \int_0^1f^2(s)\dd s$, so $G(f)$
  is simply the Wiener integral of $f\in \mathcal{L}^{\infty}$ with respect to a Brownian motion. 
 Theorem \ref{T1} thus agrees with \cite[Theorem 1.2]{BeKest}.

\item[(iii)] It would also be interesting to further investigate fluctuations around the second order limit in the super-critical case $\rho>1/2$. In the setting of Simon's model, Bertenghi \cite{Berten} observed that these are asymptotically normal; see also \cite{KuboTa}. 
 \end{enumerate} }
\end{remark}  
The number $C$ of copies of the sampled ball which are returned at a typical step does not play much of a role in  \eqref{E:9}, meaning that for such urn schemes with innovation,  reinforcement does not impact the first order asymptotic of the empirical measures of colors in the urn.
The effects of the reinforcement on the fluctuations can now  be analyzed by inspection of  the formula  \eqref{E:11} for the covariance. In short, 
consider first the case $C\equiv 0$ without reinforcement, where at each step one only returns in the urn the atoms of a point process $\xi$. By \eqref{E:9},
the empirical distribution of colors still converges to $\bar \mu$, and by Theorem \ref{T1}, the fluctuations are governed by a Gaussian bridge $G^{\mathrm{(br)}}_0$ with covariance
$$\E(G^{\mathrm{(br)}}_0(f)G^{\mathrm{(br)}}_0(g)) =  \int_{S} E_s(\xi(f)\xi(g))\bar \mu(\dd s)\qquad \text{when }\mu(f)=\mu(g)=0.$$
The comparison with \eqref{E:11} shows that the reinforced scheme with $C\not \equiv 0$
induces an additional white noise component in the fluctuations, which is seen from the term
$\bar \mu(fg) \int E_s(C^2)\bar \mu(\dd s)  $. Keep in mind that the reinforcement variable $C$ also impacts the eigenvalues  $\lambda_1$ and $\lambda_2$, and in particular that increasing $C$ while keeping the innovation $\xi$ unchanged also increases the ratio $\rho$. 

It has been briefly mentioned at the beginning of this introduction that urns schemes with infinite colors, or measure-valued P\'olya processes, have been considered in several recent works  \cite{BanJthac,  BanThac1, BanThac2, Jansoninfty, JanMV, MaiMarc, MaiVill},  and we shall now recast our contribution from this perspective.  In \cite{BanThac2} and \cite{MaiMarc}, the authors consider urns processes as  Markov chains  $(U_n)_{n\geq 0}$ of random measures (not necessarily point processes) on a measurable space $S$, such that  conditionally given $U_0, \ldots, U_n$, a random color $\sigma$ is sampled according to the probability measure  proportional to $U_n$ and then $U_{n+1}=U_n +  R_{\sigma}$, where $(R_s)_{s\in S}$ is a replacement kernel of positive measures on $S$. A more general version where now the  replacement kernels are random has been developed then in \cite{Jansoninfty} and  \cite{MaiVill}, and the framework there encompasses ours. Convergence of empirical distribution of colors has been established in \cite[Theorem 1]{MaiVill} or \cite[Section 6]{Jansoninfty}, or also \cite{BanJthac}, under certain additional conditions which however need not to be fulfilled in our setting. Fluctuations for measure-valued P\'olya processes have been considered very recently in  \cite{JanMV}, where results bearing the same flavor as our Theorem \ref{T2} are established. However \cite{JanMV} imposes a strong balance condition that would imply in our case that the total number of returned ball at each step is some constant, independently of the sampled color. Thus, although general urn schemes with random replacement kernels encompass the setting of the present work, and Theorems \ref{T:LLN},  \ref{T1} and \ref{T2} are  close relatives to results there, the present contribution cannot be reduced to  these works either.
We have chosen here to consider random replacement kernels  induced by point processes  only rather than more general random measures as in the works that were just mentioned above, mostly to keep the connection with P\'olya's original model as close as possible. 

The general route that we will follow to establish Theorems \ref{T:LLN} and  \ref{T1} is well paved since the work by Athreya and Karlin \cite{AthreyaKarlin} (see also \cite{Janson}) for urn schemes with finitely many colors, or more generally for finite-dimensional stochastic approximation \cite{BMP}. We shall exhibit natural martingales, investigate the asymptotic behavior of their quadratic variations and apply a version of the martingale central limit theorem. Here and there, we shall take a few shortcuts that where not fully available at the time when \cite{AthreyaKarlin} was written, such as stable convergence (see \cite[Section VIII.5c]{JS}). Another rather minor difference is that we do not use the embedding of urn schemes in a continuous time Markov branching process known as poissonization, but merely a Poisson subordination. More precisely, Theorem \ref{T1} will be deduced from Theorem \ref{T2}, which is a more precise  functional limit theorem for the continuous time version of the urn scheme with innovation, analogous to results of Janson \cite{Janson} in the setting of urns with finitely many colors. 

The plan for the rest of this text is as follows. Some further notation and preliminary observations on the replacement kernel and on Poisson subordination are presented in  Section 2.  The foundations of our analysis, namely the natural martingales, are laid down in Section 3; these yield Theorem \ref{T:LLN}. The large time behavior of quadratic variations is determined in Section 4. In Section 5, we establish a functional limit theorem (Theorem \ref{T2}), and finally derive Theorem \ref{T1} from the latter. 
  
  \section{Preliminaries}
\subsection{The replacement kernel}
 We introduce first some notation, referring to \cite[Chapter 1]{Kall} for background on the notions we shall use. We consider  some Borel space endowed with its Borel sigma-field $(S, \calS)$, and 
  write $\calN$ for the space of \textit{finite} counting measures on $S$; $\calN$ is also a Borel space.
  Recall that the replacement kernel $(P_s)_{s\in S}$ defines for every $s\in S$ 
 a probability distribution on $ \{-1,0,1,2,\ldots\} \times \calN$ which determines the join law of
 the number of copies $C$ and the point process of innovation $\xi$
 when the sampled ball has color $s$. We assume throughout the rest of this text and without further mention that the requirements \eqref{E:2}, \eqref{E:3} and \eqref{E:4} 
  are fulfilled.

For every $s\in S$,  $k\geq -1$ and $\Phi: \calN\to \R_+$ measurable, there is the identity
 $$E_s\left( \Phi(\xi)\ind{C=k}
 \right) = \int_{\calN} \Phi(\nu) P_s(k, \dd \nu).$$
 In particular, the intensity measure of $\xi$ under $P_s$ is given for every $f\in \mathcal{L}^{\infty}$ by 
 \begin{equation}\label{E:15}
E_s(\xi(f))= a(s)\mu(f) = \sum_{k\geq - 1} \int_{\calN} \nu(f) P_s(k,\dd \nu).
  \end{equation}
 Observe also from the notation \eqref{E:5} that
 \begin{equation}\label{E:16}
 \lambda_2= E_s(C)= \sum_{k\geq -1} k P_s(k, \calN) \qquad \text{for any }s\in S,
  \end{equation}
 and 
 \begin{align*}
 \lambda_1
 %&= \int_{S} E_s(C+\xi(a))\bar \mu(\dd s) \nonumber \\
 &= \int_{S} \bar \mu(\dd s) \left( \sum_{k\geq - 1} \int_{\calN} (k+\nu(a)) P_s(k,\dd \nu)\right). 
 \end{align*}

 We formalize the urn scheme with innovation depicted in the introduction and construct a sequence of point processes $(U_n)_{n\geq 0}$ on $S$ 
 as follows. We start from some arbitrary  $U_0\neq 0$ in $\calN$.  Next for every $n\geq 0$,   we sample a color $\sigma_n$ according to the current empirical distribution $\bar U_n$ and independently
 of $U_0, \ldots, U_{n-1}$, that is
 $$\P(\sigma_n=s\mid U_0, \ldots, U_n)=\bar U_n(\ind{s}) = U_n(\ind{s})/U_n(\mathbf{1}),\qquad s\in S.$$
 Then given $\sigma_n=s$, we  set
 \begin{equation} \label{E:17}
 U_{n+1} = U_n + C_{n+1}  \delta_{s} + \xi_{n+1},
 \end{equation}
 where  $(C_{n+1},\xi_{n+1})$ is distributed according to $P_s$ and independent of the preceding steps.
 The point process $U_n$  gives the composition of the urn after $n$ steps. 
 The sequence $(U_n)_{n\geq 0}$  is a Markov chain on 
 $\calN$, with one-step transitions given by \eqref{E:17}.

 \subsection{Poisson subordination} 
One of the most efficient and best know tools in the study of urn schemes is the so-called poissonization, which, roughly speaking, consists of
embedding the urn scheme in a continuous time multitype branching process. See for instance,  \cite{AthreyaKarlin},  \cite[Section V.9]{AthreyaNey},  \cite[Section 4.6]{Mahmoud} or \cite{Janson}. Here, we shall rather use subordination\footnote{Actually, subordination and poissonization are related one to the other by the so-called Lamperti transformation, as it will be briefly discussed in Remark \ref{R2} below. In this respect, subordination should be viewed as a mild variation of poissonization rather than a new technique in its own right. This enables us to circumvent some technical issues, such as the construction of the non-exploding branching process with general type space.} based on an independent Poisson process. The advantage is that this both circumvents the possibly delicate issue of depoissonization, and still makes some computations (notably for quadratic variations) much simpler than in discrete time. The drawback is that the independence between the counting processes of balls of different colors induced by poissonization is lost here; however this will not be an issue for us.  

Specifically, let $N=(N(t))_{t\geq 0}$ denote a standard Poisson process, so the Stieltjes measure $\dd N(t)$ is a Poisson point process on $[0,\infty)$ with intensity the Lebesgue measure.   We assume that $N$ is independent of the urn scheme
and use the jump times of $N$ as 
the set of times at which a ball is sampled from the urn scheme. Specifically, we set 
$$X_t= U_{N(t)} , \qquad t\geq 0,$$
and think of  $X=(X_t)_{t\geq 0}$ 
 as  the process describing the composition of the urn as time passes.
 We write respectively $\P$  and $\E$ for the probability measure and the mathematical expectation induced by this process, and also $\left(\mathcal F_t\right)_{t\geq 0}$ for its natural filtration.
 Occasionally, we may prefer to write $\P_m$  and $\E_m$  to stress that the initial composition of the urn is $U_0=X_0=m$ for some non-zero $m\in \mathcal N_S$. 
 
 We shall need some elementary bounds on the growth of the number of balls in the urn.
 \eject
 \begin{lemma}\label{L:invball} The following assertions hold:
 \begin{itemize}
\item[(i)]There exists some finite constant $\gamma<\infty$ such that
 $$\limsup_{t\to \infty} t^{-1} X_t(\mathbf{1}) \leq \gamma \qquad \text{a.s.}$$
 \item[(ii)] 
There exists $\varepsilon>0$ such that
 $$\lim_{t\to \infty} \e^{\varepsilon t} \P(X_t \leq \varepsilon t)= 0.$$
 \end{itemize}
 \end{lemma}
 
 \begin{proof} (i) Observe from \eqref{E:4} and Markov's inequality that 
 $$\sup_{s\in S} P_s(C+\xi(\mathbf{1}) > b) = O(b^{-2}) \qquad \text{as }b\to \infty.$$
 We infere that there is some integrable variable, whose  expectation is denoted by $\gamma$,  that dominates stochastically the number $C+\xi(\mathbf{1})$ of balls which are added to the urn at any step.
 It follows from the law of large numbers that the total number of balls in the urn growth at most linearly, 
$$
 \sup_{n\geq 1} n^{-1} U_n(\mathbf{1}) \leq \gamma \qquad\text{a.s.}
$$
 and we conclude applying the law of large numbers for the Poisson process $N$.
 
 (ii) Indeed, the process  $(X_t(\mathbf 1)-X_0(\mathbf 1))_{t\geq 0}$ which counts the number of added balls can be bounded from below by some Poisson process with rate 
 $ \inf_{s\in S} P_s(C+\xi(\mathbf{1})\geq 1)$, which is strictly positive by \eqref{E:4}. The assertion is plain from elementary large deviations estimates.
 \end{proof}

Subordination turns Markov chains into Markov jump processes. In our setting,
we get from  \eqref{E:17} that the infinitesimal generator $\mathcal G$ of  $X$ is 
given for  any bounded measurable functional  $\Phi: \calN\to \R$  and any non-zero measure  $m\in \calN$ by 
\begin{align}\label{E:18}
\mathcal G \Phi(m) &= \lim_{t\to 0+} t^{-1} \E_m(\Phi(X_t)-\Phi(X_0)) \nonumber \\
&= \int_S \bar m(\dd s) \sum_{k=-1}^{\infty} \int_{\calN} \left(\Phi(m+k\delta_s +\nu)-\Phi(m)\right) P_s(k,\dd \nu),
\end{align}
with the notation $\bar m(\dd s) = m(\dd s)/m(\mathbf{1})$. 

The starting point of our analysis based on the following observation. We shall consider linear functionals on $\calN$, that is of the type $\Phi(m)=m(f)$ for some $f\in \mathcal{L}^{\infty}$. Of course, such functionals are  not bounded as soon as the support of $f$ is infinite; nonetheless  \eqref{E:18}
still makes sense  and we have for $m\neq 0$
$$\mathcal G \Phi(m) =  \bar m\left( \mathcal Rf\right) 
,$$
where $\mathcal R: \mathcal{L}^{\infty}\to \mathcal{L}^{\infty}$ is the average replacement operator defined by \eqref{E:6}. 

In the next section, we will use these observations to exhibit simple martingales.
Recall from \eqref{E:7} and \eqref{E:8} that 
$\lambda_1$ and $\lambda_2$ appear in this setting as the two eigenvalues of $\mathcal R$.

 \section{Key martingales}

Recall that  the function $a\in {\mathcal L}^{\infty}$ appearing in \eqref{E:3} modulates the intensity of innovations. 
We set 
$$M_t(a)\coloneqq X_{t}(a)\exp\left(-\lambda_1 \int_0^{t} \frac{\dd r}{X_r(\mathbf{1})}\right), \qquad t\geq 0.$$
We also introduce for every  $f\in \mathcal{L}^{\infty}$ with  $\mu(f)=0$
$$M_t(f)\coloneqq X_{t}(f)\exp\left(-\lambda_2 \int_0^{t} \frac{\dd r}{X_r(\mathbf{1})}\right), \qquad t\geq 0.$$
Beware that the factor in the exponential is $\lambda_1$ in the first definition and $\lambda_2$ in the second.

\begin{lemma} \label{L1}   The processes $M(a)$ and $M(f)$ are square integrable martingales. 
\end{lemma}
\begin{proof} 
Recall first from \eqref{E:4} and \eqref{E:17} that
$$\E\left(|U_{n+1}(\mathbf{1})- U_{n}(\mathbf{1})|^2\right) \leq \sup_{s\in S} E_s\left(\left( C+\xi(\mathbf{1})\right)^2\right) < \infty.$$
It follows that $\E(X_t(\mathbf{1})^2)<\infty$ for every $t\geq 0$, and therefore processes $M(a)$ and $M(f)$ are square integrable.

We will henceforth focus on $M(a)$, as the proof of the martingale property for $M(f)$ relies on similar arguments. 
Our main task is to check that 
\begin{equation} \label{E:19} \text{the process }
\left(X_{t}(a) -  \lambda_1 \int_0^{t}  \bar X_r(a)  \dd r \right)_{t\geq 0}
\text{ is a martingale,}
\end{equation}
since then it follows from the It\^o formula that $M(a)$ is a martingale too.
In this direction, we   consider the bounded functional  on $\calN$ given by 
$\Phi_b(m) = m(a)\ind{|m(a)|\leq b}$ for some $b>0$. We know that 
$$\left(\Phi_b(X_t) - \int_0^t \mathcal G \Phi_b(X_r)\dd r\right)_{t\geq 0} \text{ is a martingale,}
$$
and we analyze the limit as $b\to \infty$. 

Plainly,   $\Phi_b(X_t)$ increases to $X_t(a) \in L^1(\P)$ as $b$ increases to $\infty$, and therefore
$\Phi_b(X_t)$ converges to $X_t(a)$ in $L^1(\P)$. 
Furthermore, for every $m\in \calN$, 
$$|\Phi_b(m+k\delta_s + \nu)-\Phi_b(m)| \leq\| a\|_{\infty} ( |k| + \nu(\mathbf{1})),$$
and since, thanks to \eqref{E:4},
$$\int_S \bar m(\dd s) \sum_{k=-1}^{\infty} \int_{\calN}(|k|+ \nu(\mathbf{1})) P_s(k,\dd \nu) \leq  \sup_{s\in S} E_s\left(|C| + \xi(\mathbf{1}) \right) <\infty,$$
  we deduce  from \eqref{E:18} by dominated convergence that for every $r\geq 0$,  there is the almost-sure convergence
\begin{align*}\lim_{b\to \infty} \mathcal G \Phi_b(X_r)&= 
 \int_S \bar X_r(\dd s) \sum_{k=-1}^{\infty} \int_{\calN}(ka(s) + \nu(a)) P_s(k,\dd \nu)\\
 &= \bar X_r(a) \left( \lambda_2 + \mu(a)\right) \\
 &= \lambda_1 \bar X_r(a), 
\end{align*}
where for the second identity, we used \eqref{E:15} and \eqref{E:16}.

On the other hand,  the bounds above show that for any $m\in \mathcal{N}_S$,
$$| \mathcal G \Phi_b(m)| \leq  \|a \|_{\infty}   \sup_{s\in S} E_s\left(|C| + \xi(\mathbf{1}) \right).$$
This enables us to apply dominated convergence, and for every $t\geq 0$
$$
 \lim_{b\to \infty} \int_0^t \mathcal G \Phi_b(X_r)\dd r  = \lambda_1 \int_0^t  \bar X_r(a) \dd r  \qquad \text{in }L^1(\P).
$$
This completes the proof of \eqref{E:19}. \end{proof} 

The martingale $M(a)$ is positive, hence has a
terminal value
$$M_{\infty}(a)\coloneqq \lim_{t\to \infty} M_t(a).$$

\begin{lemma} \label{L2} The martingale $M(a)$ is a  bounded in $L^2(\P)$ and $M_{\infty}(a)>0$, $\P$-a.s.
\end{lemma}
\begin{proof}  We will check that the terminal value of the oblique bracket
$\langle M(a)\rangle_{\infty}$ is integrable. 
To start with, observe from the urn scheme that for every $n\geq 0$,
\begin{align*}
\E\left(|U_{n+1}(a)-U_n(a)|^2\mid U_0, \ldots, U_n\right) &= \int_S \bar U_n(\dd s) E_s(|C a(s)+ \xi(a)|^2)\\
&\leq  \|a\|^2_{\infty} \sup_{s\in S} E_s((|C|+ \xi(\mathbf{1}))^2).
\end{align*}
Recall the assumption \eqref{E:4} and write $\beta<\infty$ for the bound above.
Since $X(a)$ is obtained from $U(a)$ by subordination via an  independent Poisson process,  the Stieltjes measure of the oblique bracket of the semimartingale  $X(a)$ satisfies
$\dd \langle X(a)\rangle_t \leq \beta\dd t$,
and then, by stochastic calculus, 
$$\langle M(a)\rangle_t \leq \beta \int_0^t   \exp\left(-2\lambda_1 \int_0^{r} \frac{\dd u}{X_u(\mathbf{1})}  \right) \dd r.$$

We next introduce the time substitution
\begin{equation} \label{E:20}
T(t)\coloneqq \inf\left \{r\geq 0:  \int_0^{r} \frac{\dd u}{X_u(\mathbf{1})}  >t\ \right\}, \qquad t\geq 0
\end{equation}
and stress that $T(t)<\infty$ a.s. for any $t>0$, thanks to  Lemma \ref{L:invball}(i). 
Write 
$$\beta^{-1} \langle M(a)\rangle_{\infty} \leq
\int_0^\infty  \exp\left(-2\lambda_1 \int_0^{r} \frac{\dd u}{X_u(\mathbf{1})}  \right) \dd r 
= \int_0^\infty    X_{T(t)}(\mathbf{1}) \e^{-2\lambda_1t} \dd t.
$$
We shall establish below that for any initial value $m\in \mathcal{N}_S$, 
\begin{equation} \label{E:21}
 \E_m\left( X_{T(t)}(\mathbf{1})\right) \leq  \frac{m(a)}{\bar \mu(a)}\e^{\lambda_1 t}.
 \end{equation}
This entails the integrability of $\langle M(a)\rangle_{\infty}$, hence the boundedness of $M(a)$ in $L^2(\P_m)$, and more precisely
\begin{equation}\label{E:22}
\E_m\left(M_{\infty}(a)^2\right) \leq m(a)^2 +  \frac{ \beta m(a)}{\bar \mu(a) \lambda_1}.
\end{equation}

Note first that since $T(t)$ is an $\mathcal F_t$-stopping time and $M(a)$ a nonnegative martingale,  $\E_m(M_{T(t)}(a))\leq m(a)$, and hence
$$ \E_m\left( X_{T(t)}(a)\right) \leq m(a) \e^{\lambda_1 t}.$$
Next, consider the function $f=a-\bar\mu(a)\mathbf{1}$, which has $\mu(f)=0$. Again by optional sampling 
applied now to the martingale $M(f)$, we have for every bounded stopping time $\tau$ that
\begin{align*}
&\bar \mu(a) \E_m\left( X_{\tau}(\mathbf{1})\exp\left( -\lambda_2\int_0^{\tau} \frac{\dd r}{X_r(\mathbf{1})}\right)\right) \\
& = \E_m\left( X_{\tau}(a)\exp\left( -\lambda_2\int_0^{\tau} \frac{\dd r}{X_r(\mathbf{1})}\right)\right) . 
\end{align*}
Apply this to $\tau=T(t)\wedge b$ and let $b\to \infty$. Then
$X_{T(t)\wedge b}(a)$ increases to $X_{T(t)}(a) \in L^1(\P_m)$, and we get by dominated convergence that
\begin{align*}
 \lim_{b\to \infty}  \E_m\left( X_{T(t)\wedge b}(a) \exp\left( -\lambda_2
\int_0^{T(t)\wedge b} \frac{\dd r}{X_r(\mathbf{1})} \right) \right)
 &= \E_m\left( X_{T(t)}(a) \e^{-\lambda_2 t}
 \right)\\
&\leq  m(a) \e^{(\lambda_1-\lambda_2) t}.
\end{align*}
By Fatou's Lemma, we get
\begin{align*}\label{E:21}
\e^{-\lambda_2 t} \E_m\left( X_{T(t)}(\mathbf{1})  \right)&\leq 
 \liminf_{b\to \infty}  \E_m\left( X_{T(t)\wedge b}(\mathbf{1}) \exp\left( -\lambda_2
\int_0^{T(t)\wedge b} \frac{\dd r}{X_r(\mathbf{1})} \right) \right)  \\
&\leq  \frac{m(a)}{\bar \mu(a)} \e^{(\lambda_1-\lambda_2) t}.
\end{align*}
The proof of \eqref{E:21} is now complete, and we turn our attention to the strict positivity of the terminal value. 

An application of Chebychev's inequality in \eqref{E:22} yields 
$$\P_m(M_{\infty}(a) = 0) \leq  \frac{ \beta }{\bar \mu(a) \lambda_1 m(a)}.$$
Observe that for any $t_0>0$, there is the equivalence 
$$M_{\infty}(a)=0 \ \Longleftrightarrow \ \lim_{t\to \infty} X_{t+t_0}(a)\exp\left( -\lambda_1 \int_{t_0}^{t_0+t} \frac{\dd r}{X_r(\mathbf{1})}\right)=0.
$$
So by the Markov property at time $t_0$,
$$\P_m(M_{\infty}(a) = 0) \leq  \frac{ \beta }{\bar \mu(a) \lambda_1 }\E_m(1/X_{t_0}(a)).$$
To conclude that $\P_m(M_{\infty}(a) = 0) =0$, we only need to let 
$t_0\to \infty$ and note that $\lim_{t_0\to \infty} X_{t_0}(a)=\infty$ a.s., since $a$ is bounded away from $0$ and the total number of balls in the urn grows to $\infty$ 
as time goes to $\infty$ a.s.
\end{proof}

We next turn our attention to the predictable quadratic variation of  the martingale $M(f)$. 

\begin{lemma} \label{L3}  Let $f\in \mathcal{L}^{\infty}$ with  $\mu(f)=0$.  The  oblique bracket of the square-integrable martingale $M(f)$ is absolutely continuous with
$$
\dd \langle M(f)\rangle_t 
=Q(X_t,f) \exp\left(-2\lambda_2 \int_0^{t} \frac{\dd r}{X_r(\mathbf{1})}\right) {\dd t},
$$
where for every non-zero $m\in \mathcal{N}_S$,
$$Q(m,f)\coloneqq  \int_{S} \bar m(\dd s)E_s\left( ( f(s) C + \xi(f))^2\right). $$
\end{lemma}

\begin{proof} 
The urn scheme \eqref{E:17} yields, in the notation of the statement, that for any  $f\in \mathcal{L}^{\infty}$ 
$$\E\left((U_{n+1}(f)-U_n(f))^2 \mid U_0,\ldots, U_n \right) = Q(U_n,f).$$
Since $X$ is obtained from $U$ by subordination via an independent Poisson process,  its optional quadratic variation process
$$[X(f)]_t\coloneqq \sum_{0< r \leq t} |\Delta X_r(f)|^2,\qquad t\geq 0,$$
has an absolutely continuous  predictable compensator $\langle X(f)\rangle$  given by 
$$\dd \langle X(f)\rangle_t = Q(X_t,f)\dd t.$$
 See \cite[\textsection I.3b]{JS} for background.

Next
 $M(f)$ is a pure jump martingale with quadratic variation
$$\dd [M(f)]_t =  \exp\left(-2\lambda_2 \int_0^{t} \frac{\dd r}{X_r(\mathbf{1})}\right)\dd [X(f)]_t .$$
The exponential term in the right-hand side is a continuous adapted process, so the
  predictable compensator $\langle M(f)\rangle$ of $[M(f)]$ is given by 
\begin{align*}\dd \langle M(f)\rangle_t &=  \exp\left(-2\lambda_2 \int_0^{t} \frac{\dd r}{X_r(\mathbf{1})}\right)\dd \langle X(f)\rangle_t\\
&= Q(X_t,f) \exp\left(-2\lambda_2 \int_0^{t} \frac{\dd r}{X_r(\mathbf{1})}\right)\dd t,
\end{align*}
see \cite[Theorem I.3.17]{JS}. This is the  formula in the statement. 
\end{proof}

Lemma \ref{L3} enables us  to bound the square moments of $M(f)$, and yields the following.

\begin{corollary}\label{C1} For every $f\in \mathcal{L}^{\infty}$, we have 
$$\lim_{t\to \infty} X_t(f) \exp\left(-\lambda_1 \int_0^{t} \frac{\dd r}{X_r(\mathbf{1})}\right) = \frac{\mu(f)}{\mu(a)} M_{\infty}(a), \qquad \text{a.s. and in }L^2(\P).$$
\end{corollary}
\begin{proof} In the case $f=a$, the claim amounts to the convergence of the martingale $M(a)$ to its terminal value; see Lemma \ref{L2}. By linearity, it thus suffices to check that for every $f\in \mathcal{L}^{\infty}$ with $\mu(f)=0$,
one has
$$\lim_{t\to \infty} M_t(f) \exp\left(-(\lambda_1-\lambda_2) \int_0^{t} \frac{\dd r}{X_r(\mathbf{1})}\right)= 0
\qquad \text{a.s. and in }L^2(\P).$$
In this direction, recall the time-substitution \eqref{E:20} and rewrite the above in the simpler form
\begin{equation} \label{E:23}
\lim_{t\to \infty}  \e^{-(\lambda_1-\lambda_2)t}M_{T(t)}(f) =0 \qquad \text{a.s. and in }L^2(\P).
\end{equation}

The assumption \eqref{E:4}
ensures that 
\begin{equation}\label{E:24}
q(f)\coloneqq \sup_{m\in \mathcal{N}_S} Q(m,f)<\infty,
\end{equation}
and we deduce from Lemma \ref{L3} 
that the time-changed process $(M_{T(t)}(f))_{t\geq 0}$ is still a square integrable martingale
with oblique bracket 
\begin{align*}\langle M(f)\rangle_{T(t)}&\leq q(f) \int_0^{T(t)}\exp\left(-2\lambda_2 \int_0^{v} \frac{\dd r}{X_r(\mathbf{1})}\right)\dd v  \\
&\leq  q(f) \int_0^t \e^{-2\lambda_2 r} {X_{T(r)}(\mathbf{1})} \dd r.
\end{align*}
The process $\e^{-(\lambda_1-\lambda_2)t}M_{T(t)}(f)$ is a supermartingale, and it therefore suffices to verify that the limit in \eqref{E:23} takes place in $L^2(\P)$. We compute
\begin{align*}\E(M_{T(t)}(f)^2) &=\E(X_0(f)^2) + 
\E(\langle M(f)\rangle_{T(t)})  \\
&\leq \E(X_0(f)^2) +  q(f) \int_0^t \e^{-2\lambda_2 r} \E\left( X_{T(r)}(\mathbf{1}) \right) \dd r.
\end{align*}
We now see from \eqref{E:21} and the fact that $\lambda_1>\lambda_2$ that 
$$\lim_{t\to \infty}  \e^{-2(\lambda_1-\lambda_2)t}\E\left( M_{T(t)}(f)^2\right) =0,$$
which completes the proof.
\end{proof}

\begin{remark}\label{R2}
\normalfont{
The time-substitution $T(t)$ used in the proofs of Lemma \ref{L2} and Corollary \ref{C1}
 is a version of the  Lamperti transformation; see \cite[Chapter 12]{Kyprianou}. It turns the measure-valued process $X$ into a continuous-time branching process with types in $S$.
In this setting of branching processes, the process $M_{T(\cdot)}(a)$ should be viewed as a Malthusian martingale, 
and Lemma \ref{L2} and Corollary \ref{C1} are then well-known;  see \cite[Theorems 7.1 and 7.2]{AthreyaNey}. }
\end{remark}

\begin{corollary}\label{C2} For every $f\in \mathcal{L}^{\infty}$, we have 
$$\lim_{t\to \infty} t^{-1} X_t(f)  = \lambda_1 \bar \mu(f), \qquad \text{a.s.}$$
\end{corollary}
\begin{proof} From Corollary \ref{C1}, since the terminal value $M_{\infty}(a)$ is almost surely finite and strictly positive, it suffices to check that
\begin{equation} \label{E:25}
\lim_{t\to \infty} t^{-1} \exp\left(-\lambda_1 \int_0^{t} \frac{\dd r}{X_r(\mathbf{1})}\right) =\frac{M_{\infty}(a)}{\lambda_1  \bar \mu(a)}, \qquad \text{a.s.}
\end{equation}
For this, we apply Corollary \ref{C1} with $f= \mathbf 1$ and get
$$\lim_{t\to \infty} \frac{1}{X_t(\mathbf 1)} \exp\left(\lambda_1 \int_0^{t} \frac{\dd r}{X_r(\mathbf{1})}\right) = \bar \mu(a)/ M_{\infty}(a) \qquad \text{a.s.}$$
This yields \eqref{E:25} by integration. 
\end{proof}

Obviously, Theorem \ref{T:LLN} follows from Corollary \ref{C2} and the law of large numbers for the Poisson process $N$.

\section{Asymptotic behaviors of quadratic variations}
 Limit theorems for martingales are of course the key
  ingredient  for Theorem \ref{T1}. Their applications require controlling quadratic variations, and for this we use mainly Lemmas \ref{L2} and \ref{L3} and \eqref{E:25}. 
 We  start with a simple observation; recall the notation \eqref{E:12}. 
 
 \begin{corollary} \label{C2'} Suppose $\rho>1/2$. Then  for every $f\in \mathcal{L}^{\infty}$ with $\mu(f)=0$, there is the convergence a.s. 
 $$\lim_{t\to \infty} t^{-\rho} X_t(f) = L'_f;$$
 and $L'_f$  is non-degenerate except if $\mu\left(|f| \right)=0$. 
 \end{corollary}
 \begin{proof} Recall the notation \eqref{E:24} and apply Lemma \ref{L3} to get
 \begin{align*}\langle M(f)\rangle_{\infty}&\leq q(f) \int_0^{\infty}\exp\left(-2\lambda_2 \int_0^t \frac{\dd r}{X_r(\mathbf{1})}\right)\dd t  \\
&\leq  q(f) \int_0^{\infty} \e^{-2\lambda_2 r} {X_{T(r)}(\mathbf{1})} \dd r.
\end{align*}
  Since $2\lambda_2>\lambda_1$, we deduce
  from \eqref{E:21} that  $\E(\langle M(f)\rangle_{\infty})<\infty$. 
 The  martingale $M(f)$ thus converges a.s. and in $L^2(\P)$. It only remains to use \eqref{E:25}. 
 \end{proof}
 
  We next turn our attention to the case $\rho\leq 1/2$ and first consider the growth of oblique brackets.  Recall the notation \eqref{E:10}.
   \begin{corollary} \label{C3} The following assertions hold  for every $f \in \mathcal L^{\infty}$ with $\mu(f)=0$:
\begin{enumerate}
\item[(i)] If  $\rho < 1/2$, then we have
$$
\lim_{t\to \infty} t^{2\rho-1} \langle M(f)\rangle_t = \frac{\sigma^2(f)}{1-2\rho } \left( \frac{M_{\infty}(a)}{  \lambda_1 \bar \mu(a)} \right)^{2\rho}\qquad \text{a.s.}
$$

\item[(ii)] If  $\rho = 1/2$, then we have
$$
\lim_{t\to \infty} \frac{1}{\log t} \langle M(f)\rangle_t =  \sigma^2(f)  \frac{M_{\infty}(a)}{ \lambda_1 \bar \mu(a)}  \qquad \text{a.s.}
$$

\end{enumerate}     
\end{corollary}
    \begin{proof}   Recall the notation of Lemma \ref{L3} and observe from Corollary \ref{C2} and the requirement $\mu(f)=0$ that 
 $$\lim_{t\to \infty} Q(X_t,f)= \sigma^2(f) \qquad \text{a.s.}$$
 We deduce from Lemma \ref{L3} and \eqref{E:25} that almost surely
 $$\frac{\dd \langle M(f)\rangle_t }{\dd t} \sim \sigma^2(f)  \left(  \frac{M_{\infty}(a)}{ t \lambda_1 \bar \mu(a)}  \right) ^{2\rho}\qquad \text{as }t\to \infty.$$
 By dominated convergence, this yields our claims. 
   \end{proof}
   
   We shall also need the following bounds on the $p$-variation of $M(f)$.
   \begin{lemma}\label{L6} Assume   \eqref{E:14}. If $\rho<1/2$, then we have for any $p\in\left(2,\alpha\wedge \frac{1}{\rho}\right)$ that
   $$
 \E \left(\sum_{r\leq t} |\Delta M_r(f)|^{p}\right)=  O( t^{1-p\rho}) \qquad \text{ as }t\to \infty.
$$
If $\rho=1/2$, then  there is some $p>2$ such that
$$
 \E \left(\sum_{t\geq 0} |\Delta M_t(f)|^p\right)< \infty .
$$
   \end{lemma}
   \begin{proof}
   Note first that for any $\varepsilon >0$ and $t\geq 0$, there is the inequality
$$ \exp\left(-2\lambda_1 \int_0^{t} \frac{\dd u}{X_u(\mathbf 1)}\right)
 \leq  \left(\frac{X_t(\mathbf 1)}{\varepsilon t}\right)^2 \exp\left(-2\lambda_1 \int_0^{t} \frac{\dd u}{X_u(\mathbf 1)}\right)+ 
 \ind{X_t(\mathbf 1)<\varepsilon t}.$$
Choosing $\varepsilon$ sufficiently small and using Lemma \ref{L:invball} and Corollary \ref{C1} for $f=\mathbf{1}$, we get
\begin{equation} \label{E:26}
\E\left ( \exp\left(-2\lambda_1 \int_0^{t} \frac{\dd u}{X_u(\mathbf 1)}\right) \right) = O(t^{-2}).
\end{equation}

Next,  computing predictable compensation in the same way as we did in the proof of Lemma \ref{L3},
we get  for any $p>0$
\begin{align*} &\E \left(\sum_{r\leq t} |\Delta M_r(f)|^p\right) \nonumber \\
&\leq \|f\|_{\infty}^p \sup_{s\in S} E_s((|C|+\xi(\mathbf 1 ))^p) \int_0^t \E\left( \exp\left(-p\lambda_2 \int_0^{r} \frac{\dd u}{X_u(\mathbf 1 )}\right) \right)\dd r.
\end{align*}
Provided that $p\lambda_2\leq 2\lambda_1$, that is $p\rho\leq 2$, 
we have by Jensen's inequality and \eqref{E:26} that 
$$\E\left( \exp\left(-p\lambda_2 \int_0^{r} \frac{\dd u}{X_u(\mathbf 1 )}\right) \right)= O(t^{-p\rho}).$$

The statement is now clear when $\rho<1/2$, and in the case $\rho = 1/2$, we may just  take  any  $p\in(2, \alpha]$.    \end{proof}

 The asymptotic behavior of optional quadratic variations can be deduced in turn from Corollary \ref{C3} and the following observation.

    \begin{lemma} \label{L5} Assume   \eqref{E:14} and that $\rho \leq 1/2$. Then for every $f\in \mathcal{L}^{\infty}$ with $\mu(f)=0$, we have
    $$\lim_{t\to \infty}\frac{ [M(f)]_t}{ \langle M(f)\rangle_t} = 1 \qquad \text{in probability.}$$
    \end{lemma}
\begin{proof} Suppose first that $\rho<1/2$. 
By Corollary \ref{C3},  it suffices to verify that
\begin{equation}\label{E:27}
\lim_{t\to \infty} t^{2\rho-1} \E( | D_t(f)| )=  0 ,
\end{equation}
where    $D(f)\coloneqq [M(f)]-\langle M(f)\rangle$. Since $\langle M(f)\rangle$ is continuous, $D$ is a pure jump martingale whose jumps 
coincide with those of $[M(f)]$, so $|\Delta D|^2 = |\Delta M(f)|^4$. We can bound the quadratic variation of $D$ for any $2< p \leq 4$ by
$$
[D]_t \leq \left(\sum_{r\leq t} |\Delta M_r(f)|^p\right)^{4/p}.
$$
Choosing $p$ as in Lemma \ref{L6},  and applying the Jensen and the Burkholder-Davis-Gundy inequalities, we get
$$\E(|D_t|) = O(t^{2(1-p\rho )/p})= o(t^{1-2\rho});$$
so \eqref{E:27} is verified. 

In the case $\rho=1/2$, the same argument shows that the martingale $D$ is bounded in $L^{p/2}(\P)$ and hence converges a.s.
\end{proof}

 \section{Applying stable limit theorems}
 
The purpose of this section is to establish the following functional limit theorem from which Theorem \ref{T1}(ii-iii) will  readily follow. 
The space $\mathbb D(\R_+,\R)$ of c\`adl\`ag real-valued paths is endowed with the Skorokhod topology; see \cite[Section VI.1]{JS}. 
Recall the notation \eqref{E:10}, and let also $W=(W(t))_{t\geq 0}$ denote a standard real Brownian motion.

 \begin{theorem}\label{T2} Assume  \eqref{E:14} and take any  $f \in \mathcal L^{\infty}$ with $\mu(f)=0$.
 
 \begin{enumerate} 
 \item[(i)] If $\rho<1/2$, then 
 there is the convergence in distribution on $\mathbb D(\R_+,\R)$ in the sense of Skorokhod
 $$\lim_{n\to \infty} \left(\frac{X_{tn}(f)}{ \sqrt{n}} \right)_{t\geq 0} = \left(   \frac{\sigma (f)}{\sqrt{1-2\rho}}\, t^{\rho} W\left( t^{1-2\rho} \right) \right)_{t\geq 0}.$$

\item[(ii)] If $\rho=1/2$, then 
 there is the convergence in distribution on $\mathbb D(\R_+,\R)$ in the sense of Skorokhod
  $$\lim_{n\to \infty} \left(\frac{X_{\e^{nt}}(f)}{ \sqrt{n \e^{nt}}} \right)_{t\geq 0} = \left( \sigma (f) W(t) \right)_{t\geq 0}.$$
 \end{enumerate} 
 \end{theorem}
  
 \begin{proof}[Proof of Theorem \ref{T2}]  Again, we shall only prove (i) as the argument in the critical case $\rho=1/2$ requires only simple modifications.  The claim is plain from Corollary \ref{C3}(iii) when $\sigma^2(f)=0$, so we shall also assume henceforth $\sigma^2(f)\neq 0$. 
The proof relies on general results on the so-called stable convergence for martingales (see \cite[Section VIII.5.c]{JS} for background) and the technical results we established so far.  
 
 We shall apply \cite[Theorem  VIII.5.50]{JS} 
 for  the martingale $Y$ resulting from $M(f)$ after the deterministic time change $t\mapsto t^{1/(1-2\rho)}$, namely
 $$Y_t\coloneqq M_{t^{1/(1-2\rho)}}(f), \qquad t\geq 0.$$
 More precisely, the statement there requests, for the sake of simplicity, boundedness of the jumps of $Y$ (see \cite[Eq. VIII.5.48]{JS}); however its proof only uses the weaker condition VIII.3.14, which, in the present setting, follows from Lemma \ref{L6}. 
Next, we have from Corollary \ref{C3}(i) and Lemma \ref{L5},
$$\lim_{n\to \infty} n^{-1} [Y]_{nt} =\lim_{n\to \infty} n^{-1} [M]_{(nt)^{1/(1-2\rho)}} = \eta^2 t\qquad \text{in probability,}$$
where 
 $$\eta= 
    \frac{ \sigma(f)}{\sqrt{1-2\rho}}  \left(\frac{M_{\infty}(a)}{\lambda_1 \bar \mu(a)}\right)^{\rho}. $$
 This is  the condition VIII.5.49 in \cite{JS}.  
        
  Theorem VIII.5.50 and Proposition VIII.5.33(ii) from \cite{JS} now entail that the sequence of pairs $\left(\eta,(n^{-1/2} Y_{nt})_{t\geq 0}\right)$ converges in law to $(\eta,\eta W)$,
  where $W$ is a standard Brownian motion independent of $\eta$. Recall also from Lemma \ref{L2}  that $\eta>0$ a.s. 
  It follows that 
  $$\frac{ n^{-(1-2\rho)/2} }{\eta} Y_{(nt)^{1-2\rho}}= \frac{  X_{nt}(f) }{ \sqrt n}\times   \frac{ n^{\rho} }{\eta} 
  \exp\left(-\lambda_2 \int_0^{nt} \frac{\dd r}{X_r(\mathbf{1})}\right), \quad t\geq 0,$$ converges in law to the time-changed Brownian motion
  $\left(W(t^{1-2\rho})\right)_{t\geq 0}$. 
  To complete the proof, it suffices to recall from \eqref{E:25}   that
  $$\lim_{n\to \infty}  \frac{ n^{\rho} }{\eta} 
  \exp\left(-\lambda_2 \int_0^{nt} \frac{\dd r}{X_r(\mathbf{1})}\right)
  = \frac{1}{\eta} \left( \frac{M_{\infty}(a)}{t \lambda_1 \bar \mu(a)} \right)^{\rho} =  \frac{\sqrt{1-2\rho}}{ \sigma(f)} t^{-\rho},
  $$
  and apply Slutsky's lemma.
  \end{proof}

Theorem \ref{T1} now follows readily from Corollary \ref{C2} and  Theorems  \ref{T2} and \ref{T:LLN}. 
 \begin{proof}[Proof of Theorem \ref{T1}]
 (i) By linearity, it suffices to consider the case when $\mu(f)=0$. We 
 write first
 $$n^{1-\rho} \bar U_n(f) = \frac{n}{U_n(\mathbf{1})} \times n^{-\rho} U_n(f),$$
 and recall from  Theorem \ref{T:LLN} that the first term in the product in the right-hand side converges to $1/\lambda_1$ as $n\to \infty$. 
 We next rewrite  Corollary \ref{C2}  in the form
  $$\lim_{t\to \infty} \left(\frac{N(t)}{t}\right)^{\rho}N(t)^{-\rho} U_{N(t)}(f) = L'_f,\qquad \text{a.s.}$$
 Since $N(t)\sim t$ a.s., this yields $n^{-\rho} U_n(f)\sim L'_f$ and the claim is proved with $L_f= \lambda^{-1} L'_f$. 
 
 (ii)  Fix first $f\in \mathcal L^{\infty}$ with $\mu(f)=0$. Write $0<\gamma(1) <\gamma(2)< \ldots$ for the increasing sequence of the jump times of the Poisson process $N$, 
 so $N(\gamma(n))=n$ and $U_n(f)= X_{\gamma(n)}(f)$. Since  $\lim_{n\to \infty} n^{-1}\gamma(n)=1$ a.s., we deduce from Theorem \ref{T2}(i) that  there is the convergence
 $$\lim_{n\to \infty}\left( \frac{1-2\rho}{n}\right)^{1/2} U_n(f)=G(f)=G^{\mathrm{(br)}}(f) \qquad \text{ in distribution},$$
  where $G$ is the centered Gaussian process on $\mathcal L^{\infty}$ with covariance given by \eqref{E:11}
  and $G^{\mathrm{(br)}}$ its bridge version. It now follows from Theorem \ref{T:LLN} and  Slutsky's lemma
  that 
 $$\lim_{n\to \infty}  \lambda_1 \left( (1-2\rho)n\right)^{1/2}\bar U_n(f)=G^{\mathrm{(br)}}(f) \qquad \text{ in distribution.}$$
Since $\bar U_n$ and $G^{\mathrm{(br)}}$ are both translation invariant
 (in the sense that $\bar U_n(f+c)= \bar U_n(f)$ and 
 $G^{\mathrm{(br)}}(f+c)= G^{\mathrm{(br)}}(f)$ for any $c\in\R$),
  the convergence in distribution above holds for any fixed $f\in \mathcal L^{\infty}$, without the restriction that $\mu(f)=0$. 
  Finally, that the latter holds more generally jointly for any finite family of functions in $\mathcal L^{\infty}$ stems from the Cram\'er-Wold device. 
  
 (iii) The argument is similar to  (ii). 
  \end{proof}

\bibliography{Polya.bib}

\end{document}